\documentclass{amsart}
\usepackage{amsmath}
\usepackage{amssymb}
\usepackage{mathrsfs}

\newtheorem{pro}{Proposition}

\newtheorem{cor}{Corollary}
\newtheorem*{nota}{Remark}
\newtheorem*{remarks}{Remarks}
\newcommand{\p}{\mathbb{P}}
\newcommand{\e}{\mathbb{E}}
\newcommand{\re}{\mathbb{R}}
\date{\today}

\begin{document}
\title{Radial Dunkl processes associated with dihedral systems} 
\maketitle
\centerline{NIZAR DEMNI  \footnote{Adress: Fakult\"at f\"ur mathematik, universit\"at Bielefeld. Postfach 100131, Bielefeld, Germany.  Email: demni@math.uni-bielefeld.de}}
\begin{abstract}
We give some interest in radial Dunkl processes associated with dihedral systems. We write down the semi group density and as a by-product the generalized Bessel function and the $W$-invariant generalized Hermite polynomials. Then, a skew product decomposition, involving only independent Bessel processes, is given and the tail distribution of the first hitting time of boundary of the Weyl chamber is computed.
\end{abstract} 

\section{A quick reminder}
We refer the reader to \cite{Dunkl} and \cite{Hum} for facts on root systems and to \cite{Chy}, \cite{Ros} for facts on radial Dunkl processes. Let $R$ be a reduced root system in a finite euclidean space $(V,<,>)$ with positive system $R_+$ and simple system $S$. Let $W$ be its reflection group and $C$ be its positive Weyl chamber. The radial Dunkl process $X$ associated with $R$ is a continuous paths Markov process valued in $\overline{C}$ whose generator acts on $C^2(\overline{C})$-functions as 
\begin{align*}
\mathscr{L}_k u(x) = \frac{1}{2}\Delta u(x) + \sum_{\alpha \in R_+}k(\alpha)\frac{<\nabla u(x),\alpha>}{<x,\alpha>}
\end{align*}
with $<\nabla u(x),\alpha> = 0$ whenever $<x,\alpha> =0$, where $\Delta,\nabla$ denote the euclidean Laplacian and the gradient respectively and $k$ is a positive multiplicity function, that is, a $\re_+$-valued $W$-invariant function. The semi group density of $X$ with respect to the Lebesgue measure in $V$ is given by
\begin{equation}\label{RSG}
p_t^k(x,y) =\frac{1}{c_kt^{\gamma + m/2}} e^{-(|x|^2 + |y|^2)/2t}D_k^W\left(\frac{x}{\sqrt t},\frac{y}{\sqrt t}\right) \omega_k^2(y), \quad x,y \in \overline{C}
\end{equation} 
where $\gamma = \sum_{\alpha \in R_+}k(\alpha)$ and $m= \textrm{dim} V$ is the rank of $R$. The weight function $\omega_k$ is given by 
\begin{equation*}
\omega_k(y) = \prod_{\alpha \in R_+} <\alpha,y> ^{k(\alpha)}
\end{equation*}
and $D_k^W$ is the generalized Bessel function. Thus, $\mathscr{L}_k$ may be written as
\begin{equation}\label{gen}
\mathscr{L}_k u(x) = \frac{1}{2}\Delta u(x) + <\nabla u(x), \nabla \log \omega_k(x)>.
\end{equation}

\section{Motivation}
Several reasons motivated us to investigate radial Dunkl processes associated with dihedral root systems. First, the dihedral group is a Coxeter, yet non Weyl in general, reflections group and covers an exceptional Weyl group known in the literature as $G_2$ which is of a particular interest (\cite{Baez}). Second, the study of the Dunkl operators associated with dihedral root systems revealed a close relation to Gegenbauer and Jacobi polynomials which have interesting geometrical interpretations such as being spherical harmonics and eigenfunctions of the radial part of the Laplacian on the sphere (\cite{Bakry}, \cite{Dunkl}, \cite{Dunkl1}). The latter operator is a particular case of the Jacobi operator which generates a diffusion known as the Jacobi process that may be represented, up to a random time change, by means of two independent Bessel processes (\cite{War}). Since the norm of the radial Dunkl process is a Bessel process, we wanted to gather all these materials in the present work and to see how do they interact.  The last reason is that \cite{Dem} and \cite{Dem1} emphasize the irreducible root systems of types $A,B,C,D$ which, together with dihedral root systems, exhaust the infinite families of irreducible root systems associated with finite Coxeter groups.\\
The remanining part of the paper consists of five sections. In order to be self-contained, some needed facts on dihedral systems are collected in the next section. 
Then, we write down the semi group density: to proceed, we perform a detailed analysis of the so-called spherical motion (see below). As a by-product, we deduce the generalized Bessel function. Once we did, we express the $W$-invariant counterparts of the generalized Hermite polynomials (\cite{Ros}) as a product of univariate Laguerre and Jacobi polynomials. Next, we give a skew product decomposition of the radial Dunkl process using only independent Bessel processes. This mainly follows from the skew-product decomposition of the Jacobi process derived in \cite{War}. Finally, we compute the tail distribution of the first hitting time of $\partial C$. 
\section{Dihedral groups and Dihedral systems}
The dihedral group, denoted by $\mathcal{D}_2(n)$ for $n \geq 3$, is defined as the group of orthogonal transformations that preserve a regular $n$-sided polygon in $V=\re^2$ centered at the origin. Without loss of generality, one may assume that the $y$-axis is a mirror for the polygon. It contains $n$ rotations through multiples of $2\pi/n$ and $n$ reflections about the diagonals of the polygon. By a diagonal, we mean a line joining two opposite vertices or two midpoints of opposite sides if $n$ is even, or a vertex to the midpoint of the opposite side if $n$ is odd. The corresponding dihedral root system, $I_2(n)$, is characterized by its positive and simple systems  given by: 
\begin{equation*}
R_+ = \{-i e^{i\pi l/n} := -ie^{i\theta_l},\, 1 \leq l \leq n\},\quad S = \{e^{i\pi/n}e^{-i\pi/2}, e^{i\pi/2}\}
\end{equation*}  
so that the Weyl chamber is a wedge of angle $\pi/n$. The reader can check that, for instance, $I_2(3)$ (equilateral triangle-preserving) is isomorphic to $R=A_2$ and $I_2(4)$ 
(square-preserving) is nothing but $R=B_2$ (see \cite{Hum} for the definitions of both root systems). However, it is a bit more delicate to see that $I_2(6)$ (hexagon-preserving) corresponds to the exceptional Weyl group $G_2$ 
(\cite{Baez}). When $n = 2p, p\geq 2$, there are two orbits so that $k=(k_0,k_1) \in \re_+^2$, otherwise, there is only one orbit and $k$ takes only one positive value. 

\section{Semi group density}
\subsection{Spherical motion}
Let us recall the skew-product decomposition of the radial Dunkl process into a radial and a spherical parts (\cite{Chy} p. 53):
\begin{equation*}
(X_t)_{t \geq 0} = \left(|X_t|\Theta_{A_t}\right)_{t \geq 0} \quad A_t = \int_0^t \frac{ds}{|X_s|^2},
\end{equation*} 
where $|X|$ is a Bessel process of index $\gamma$ (\cite{Rev}) and $\Theta$ is the spherical motion of $X$ valued in the intersection of the sphere of the euclidean space $V$ and the closure $\overline{C}$ of the positive Weyl Chamber, and is independent of $|X|$. For dihedral systems, $\Theta$ is valued in the unit circle and may be written as $(\cos\theta,\sin \theta)$ for some process $\theta$ valued in the interval $[0,\pi/n]$. 

\subsection{Semi group density of $\theta$} 
Let us  first split $\mathscr{L}_k$ of $X$ into a radial and a spherical parts. To proceed, we use (\ref{gen}) together with the expression of $\omega_k$ in polar coordinates (up to a constant factor, \cite{Dunkl} p. 205) 
\begin{eqnarray*}
\omega_k(r,\theta) &=& r^{nk} (\sin n\theta)^k \\
\omega_k(r,\theta) &=& r^{p(k_0+k_1)}[\sin (p\theta)]^{k_0}[\cos (p\theta)]^{k_1}
\end{eqnarray*}
for $n$ being odd, $n=2p$ respectively and $y = (r,\theta) \in \re^2$. Thus 
\begin{equation*}
\mathscr{L}_k  = \frac{1}{2}\left[\partial_r^2  + \frac{2\gamma +1}{r}\partial_r\right]  + \frac{1}{r^2}\left[\frac{\partial_{\theta}^2}{2} +nk\cot(n\theta)\partial_{\theta}\right]
\end{equation*}
when $n$ is odd, where $\gamma = nk$, and
\begin{equation*}
\mathscr{L}_k = \frac{1}{2}\left[\partial_r^2  + \frac{2\gamma +1}{r}\partial_r\right]  + \frac{1}{r^2}\left[\frac{\partial_{\theta}^2}{2} +p(k_0\cot(p\theta) - k_1\tan(p\theta))\partial_{\theta}\right]
\end{equation*}
when $n=2p$, where $\gamma = p(k_0+k_1)$. It follows that the generator of $\theta$, say $\mathscr{L}_k^{\theta}$, acts on smooth functions as
\begin{eqnarray*}
\mathscr{L}_k^{\theta} &=& \frac{\partial_{\theta}^2}{2} +nk\cot(n\theta)\partial_{\theta}\\
\mathscr{L}_k^{\theta} &=& \frac{\partial_{\theta}^2}{2} +p(k_0\cot(p\theta) - k_1\tan(p\theta))\partial_{\theta}
\end{eqnarray*}
when $n$ is odd, $n=2p$ respectively. Now, it is easy to see that the process $N$ defined by $N_t := n\theta_{t/n^2}$ satisfies 
\begin{equation*}
dN_t= dB_t + k\cot (N_t) dt 
\end{equation*}
when $n$ is odd, while $(M_t := p\theta_{t/p^2})_{t \geq 0}$ satisfies 
\begin{equation*}
dM_t = dB_t  + [k_0\cot(M_t) - k_1\tan(M_t)]dt
\end{equation*}
when $n=2p$, $B$ being a real Brownian motion. Let us first investigate the case of even $n=2p$.
The generator of $M$ has a discrete spectrum given by $\lambda_j = -2j(j+k_0+k_1), j \geq 0$ corresponding to the Jacobi polynomials 
$P_j^{k_1 -1/2, k_0-1/2}(\cos(2\theta))$ (see \cite{Dunkl}, p. 201). It is known that this set of orthogonal eigenpolynomials is complete for the Hilbert space 
$L^2([0,\pi/2], \mu_k(\theta)d\theta)$ where 
\begin{equation*}
\mu_k(\theta) := c_{k} \sin(\theta)^{2k_0}\cos(\theta)^{2k_1}
\end{equation*}
for some constant $c_{k}$. Accordingly, $M$ has a semi group density, say $m_t^k(\phi,\theta)$, given by (we use orthonormal polynomials, \cite{Schou} p. 29)
\begin{equation}\label{Jacobi}
m_t^k(\phi,\theta) = \sum_{j\geq 0}e^{\lambda_j t}P_j^{k_1 -1/2, k_0-1/2}(\cos(2\phi))P_j^{k_1 -1/2, k_0-1/2}(\cos(2\theta))\mu_k(\theta)
\end{equation}
for $\phi,\theta \in [0,\pi/2]$. It follows that the semi group density of $\theta$, say $K_t^{k,p}$, is given by 
\begin{equation*}
K_t^{k,p}(\phi,\theta) = p m_{p^2t}^k(p\phi,p\theta), \quad \phi,\theta \in [0,\pi/(2p)].
\end{equation*}
A similar spectral description holds for odd $n$: the generator of $N$ has a discrete spectrum given by $\lambda_j = -2j(j+k), j \geq 0$ corresponding to
$P_j^{-1/2, k-1/2}(\cos(2\theta))$. 
 
\subsection{Semi group density of $X$}
Let $(r,\theta) \mapsto f(r,\theta)$ be a nice function and let $\p_{\rho,\phi}$ denote the law of $X$ starting at $x = (\rho,\phi) \in C$. Then, using the independence of $\theta$ and $|X|$ together with Fubini's Theorem, one has
\begin{align*}
&\mathbb{E}_{\rho,\phi}[f(|X_t|,  \theta_{A_t})] = \mathbb{E}_{\rho,\phi}[\mathbb{E}_{\rho,\phi}[f(|X_t|,  \theta_{A_t})| \sigma(|X_s|, s \leq t)]] 
\\& = \int_0^{\pi/(2p)}\sum_{j \geq 0}\mathbb{E}_{\rho}^{\gamma}[f(|X_t|, \theta) e^{\lambda_j p^2A_t}]P_j^{l_1, l_0}(\cos(2p\phi))P_j^{l_1, l_0}(\cos(2p\theta))\mu_k(p\theta)
d\theta
\end{align*}
where $\p_{\rho}^{\gamma}$ is the law of the Bessel process $|X|$ starting at $\rho$ and of index $\gamma$. Next, for every $\theta \in [0,\pi/(2p)]$
\begin{align*}
\mathbb{E}_{\rho}^{\gamma}[f(|X_t|, \theta) e^{\lambda_j p^2A_t}] &= \mathbb{E}_{\rho}^{\gamma}[\mathbb{E}_{\rho}^{\gamma}[f(|X_t|, \theta) e^{\lambda_j p^2A_t}| |X_t|]
 \\&= \int_0^{\infty}\mathbb{E}_{\rho}^{\gamma}[e^{\lambda_j p^2A_t}| |X_t|=r] f(r,\theta)q_t(\rho,r) dr
\end{align*}
where $q_t(\rho,r)$ is the semi group density of the Bessel process $|X|$ of index $\gamma$ (\cite{Rev}):
\begin{equation*}
q_t(\rho,r) = \frac{1}{t} \left(\frac{r}{\rho}\right)^{\gamma}r e^{-(\rho^2+r^2)/2t} {\it I}_{\gamma}\left(\frac{\rho r}{t}\right)
\end{equation*}
where ${\it I}_{\gamma}$ is the modified Bessel function of index $\gamma$ (\cite{Leb} p. 108). Moreover (see \cite{Yor} p. 80)
\begin{equation*}
\mathbb{E}_{\rho}^{\gamma}[e^{\lambda_j p^2A_t}| |X_t|=r]  = \frac{{\it I}_{\sqrt{\gamma^2 - 2\lambda_j p^2}}(\rho r/t)}{{\it I}_{\gamma}(\rho r/t)},\, \lambda_j = -2j(j+k_0+k_1). 
\end{equation*}
Thus, we proved that
\begin{pro}
The semi group density of the radial Dunkl process associated with even Dihedral groups $\mathcal{D}_2(2p)$ is given by  
\begin{align*}
p_t^k(\rho,\phi,r, \theta) &= \frac{1}{c_kt} \left(\frac{r}{\rho}\right)^{\gamma} e^{-(\rho^2+r^2)/2t}\sin^{2k_0}(p\theta)\cos^{2k_1}(p\theta) 
\\& \sum_{j \geq 0}{\it I}_{2jp+\gamma}\left(\frac{\rho r}{t}\right)P_j^{l_1, l_0}(\cos(2p\phi))P_j^{l_1, l_0}(\cos(2p\theta))
\end{align*}
with respect to $dr d\theta$, where $c_k$ is a normalizing constant, $l_0 = k_0 - 1/2, l_1 = k_1-1/2$ and $\rho, r \geq 0, 0 \leq \phi,\theta \leq \pi/2p$. For odd dihedral systems, one has to substitute in the above formula $k_1=0, k_0 = k, p=n$.
\end{pro}

\begin{remarks}
1/The $j$-th Jacobi polynomial $P_j^{k_1-1/2,k_0-1/2}(\cos(2p\theta))$ can be replaced by the generalized Gegenbauer polynomial $C_{2j}^{k_1,k_0}(\cos (p\theta))$ 
(see \cite{Dunkl}, p. 27). For $k_1=0$, $C_{2j}^{k_1,k_0}(\cos (p\theta))$ reduces to the Gegenbauer polynomial $C_{2j}^{k_0}(\cos (p\theta))$.\\
2/The heat kernel corresponding to a planar Brownian motion starting at $x \in C$ and reflected on $\partial C$ corresponds to $k \equiv 0$. Using the above formula, one deduces  
\begin{align*}
p_t^0(\rho,\phi,r, \theta) &= \frac{1}{c_0t}e^{-(\rho^2+r^2)/2t} \sum_{j \geq 0}{\it I}_{2jp}\left(\frac{\rho r}{t}\right)T_j(\cos(2p\phi))T_j(\cos(2p\theta))
\end{align*}
where $T_j$ is the orthonormal $j$-th Tchebycheff polynomial defined by
\begin{equation*}
T_j(\cos \theta) = \cos(j\theta),\quad j \geq 0.
\end{equation*}
Thus 
\begin{align*}
p_t^0(\rho,\phi,r, \theta) &= \frac{1}{c_0t}e^{-(\rho^2+r^2)/2t} \sum_{j \geq 0}{\it I}_{2jp}\left(\frac{\rho r}{t}\right)\cos(2jp\phi)\cos(2jp\theta).
\end{align*}
For $k \equiv 1$, one recovers the Brownian motion conditioned to stay in a wedge of angle $\pi/n$ which is the $h = \omega_1$-transform in Doob's sense of a planar Brownian motion killed when it first hits $\partial C$ (\cite{Gra}). More precisely, one has for $n=2p$
\begin{align*}
p_t^1(\rho,\phi,r, \theta) &= \frac{ \omega_1^2(r,\theta)}{c_1t} \left(\frac{1}{r\rho}\right)^{2p} e^{-(\rho^2+r^2)/2t}
 \sum_{j \geq 0}{\it I}_{2(j+1)p}\left(\frac{\rho r}{t}\right)U_j(\cos(2p\phi))U_j(\cos(2p\theta)),
\end{align*}
where $U_j$ is the $j$-th Tchebycheff polynomial of the second kind defined by
\begin{equation*}
U_j(\cos \theta) = \frac{\sin(j+1)\theta}{\sin \theta}, j \geq 0.
\end{equation*}

\begin{equation*}
\omega_1(r,\theta) = cr^{2p}\sin(2p\theta)
\end{equation*}
for some constant $c$, elementary computations yield
\begin{align*}
p_t^1(\rho,\phi,r, \theta) &= \frac{\omega_1(r,\theta)}{\omega_1(\rho,\phi)} \frac{e^{-(r^2+\rho^2)/2t}}{c_1t}
\sum_{j \geq 0}{\it I}_{2(j+1)p}\left(\frac{\rho r}{t}\right)\sin[2p(j+1)\phi)]\sin[2(j+1)p\theta]
\end{align*}
which agrees with the $\omega_1$-transform property. Besides, one deduces that the semi group density of a planar Brownian motion killed when it first hits the $\partial C$ is given by  
\begin{equation*}
p_t^C(\rho,\phi,r,\theta) = \frac{e^{-(r^2+\rho^2)/2t}}{t}\sum_{j \geq 0}{\it I}_{2(j+1)p}\left(\frac{\rho r}{t}\right)\sin[2p(j+1)\phi)]\sin[2p(j+1)\theta].
\end{equation*}
The expression of $p_t^C$ should be compared with Lemma 1 in \cite{Ban}. 
\end{remarks} 

Writing $p_t^k$ as
\begin{align*}
p_t^k(\rho,\phi,r, \theta) &= \frac{1}{c_kt^{\gamma +1}} \left(\frac{t}{r\rho}\right)^{\gamma} e^{-(\rho^2+r^2)/2t} \omega_k^2(r,\theta)
\\& \sum_{j \geq 0}{\it I}_{(2j+k_0+k_1)p}\left(\frac{\rho r}{t}\right)P_j^{l_1, l_0}(\cos(2p\phi))P_j^{l_1, l_0}(\cos(2p\theta))
\end{align*}
and with regard to (\ref{RSG}), we are led to the by-product

\begin{cor}[Generalized Bessel function]\label{GBF}
For even dihedral groups, the generalized Bessel function is given by 
\begin{equation*}
D_k^W(\rho,\phi,r,\theta) = c_{p,k}\left(\frac{2}{r\rho}\right)^{\gamma}
\sum_{j \geq 0}{\it I}_{2jp +\gamma}(\rho r)P_j^{l_1, l_0}(\cos(2p\phi))P_j^{l_1, l_0}(\cos(2p\theta))
\end{equation*}
where $\gamma = p(k_0+k_1)$. For odd dihedral groups, one has  
\begin{equation*}
D_k^W(\rho,\phi,r,\theta) = c_{n,k}\left(\frac{2}{r\rho}\right)^{\gamma}
\sum_{j \geq 0}{\it I}_{2jn+\gamma}(\rho r)P_j^{-1/2, l_0}(\cos(2n\phi))P_j^{-1/2, l_0}(\cos(2n\theta))
\end{equation*}
where $\gamma = nk, l= k-1/2$. The constant $c_{p,k}$ and $c_{n,k}$ are such that $D_k^W(0,0,r,\theta) = |W|$.
\end{cor}

\section{$W$-invariant generalized Hermite polynomials}
In this section, we shall express the $W$-invariant counterparts of the so-called generalized Hermite polynomials by means of univariate Laguerre and Jacobi polynomials. 
This is done in three steps.

\subsection{Generalized Hermite polynomials}
Recall from \cite{Ros} that the generalized Hermite polynomials $(H_{\tau})_{\tau \in \mathbb{N}^m}$  are defined by
\begin{equation*}
H_{\tau}(x) = [e^{-\Delta_k/2} \phi_{\tau}](x)
\end{equation*}
where $\Delta_k$ denotes the Dunkl Laplacian (\cite{Ros}) and $(\phi_{\tau})_{\tau \in \mathbb{N}^m}$ is a basis of homogeneous polynomials orthogonal with respect to the pairing inner product defined in \cite{Dunkl01} (see also \cite{Ros}):
\begin{equation*}
[p,q]_k = \int_V e^{-\Delta_k/2}p(x)e^{-\Delta_k/2}q(x) \omega_k^2(x) dx
 \end{equation*}
 for two polynomials $p,q$ (up to a constant factor). The family $(H_{\tau})_{\tau}$ is then said to be associated to the basis $(\phi_{\tau})_{\tau}$. Their $W$-invariant counterparts are defined by 
 \begin{equation*}
H_{\tau}^W(x) := \sum_{w \in W}H_{\tau}(wx).
\end{equation*}

\subsection{Mehler-type formula}
It is known that $(H_{\tau})_{\tau}$ satisfies a Mehler-type formula (\cite{Dunkl} p. 246\footnote{we use a different normalization from the one used in \cite{Dunkl}.})
\begin{align*}
\sum_{\tau \in \mathbb{N}^m}H_{\tau}(x)H_{\tau}(y)r^{|\tau|} = \frac{1}{(1-r^2)^{\gamma + m/2}}\exp{-\frac{r^2(|x|^2+|y|^2)}{2(1-r^2)}} D_k\left(x,\frac{r}{1-r^2}y\right) 
\end{align*}
for $0 < r < 1$, $x,y \in V$. Analogous formula is satisfied by $(H_{\tau}^W)_{\tau}$ and follows after summing twice over $W$ and using $D_k(wx,w'y) = D_k(x,w^{-1}w'y)$ (\cite{Ros}): 
\begin{align}\label{Mehler}
\sum_{\tau \in \mathbb{N}^m}H_{\tau}^W(x)H_{\tau}^W(y)r^{|\tau|} = \frac{|W|}{(1-r^2)^{\gamma + m/2}}\exp{-\frac{r^2(|x|^2+|y|^2)}{2(1-r^2)}} D_k^W\left(x,\frac{r}{1-r^2}y\right). 
\end{align}

\subsection{Dihedral systems}
Let us express $D_k^W$ through the hypergeometric function ${}_0\mathscr{F}_1$. This is done via the relation (\cite{Leb})
\begin{equation*}
{\it I}_{\nu}(z) = \frac{1}{\Gamma(\nu+1)}z^{\nu}{}_0\mathscr{F}_1(\nu+1,z).
\end{equation*}
It follows that 
\begin{equation*}
D_k^W(\rho,\phi,r,\theta) = c_{p,k}\sum_{j \geq 0}\frac{(\rho r/2)^{2jp}}{\Gamma(2jp+\gamma+1)}{}_0\mathscr{F}_1\left(2jp + \gamma +1, \frac{\rho^2 r^2}{4}\right)
P_j^{l_1, l_0}(\cos(2p\phi))P_j^{l_1, l_0}(\cos(2p\theta)).
\end{equation*}
Using the Mehler-type formula for univariate Laguerre polynomials (\cite{Bak} p. 200):
\begin{align*}
\sum_{q \geq 0}\frac{q!}{(2jp+\gamma+1)_q}&L_q^{2jp+\gamma}(\rho^2/2)L_q^{2jp +\gamma} (r^2/2) z^{2q} = (1-z^2)^{-2jp - \gamma -1}
\\& e^{-z^2(\rho^2 + r^2)/[2(1-z^2)]}{}_0\mathscr{F}_1\left(2jp+\gamma +1, \frac{z^2\rho^2r^2}{4(1-z^2)^2}\right), \quad |z| < 1,
\end{align*}
one gets
\begin{align*}
&(1-z^2)^{-\gamma -1}e^{-z^2(\rho^2 + r^2)/[2(1-z^2)]} D_k^W\left(\rho,\phi,\frac{zr}{1-z^2},\theta\right) =  \\& c_{p,k} 
\sum_{j,q \geq 0}\frac{q!}{\Gamma(2jp+q+\gamma+1)} N_{j,p}^{k,W}(\rho,\phi) N_{j,p}^{k,W}(r,\theta)\left(\frac{\rho r}{2}\right)^{2jp}z^{2(q+jp)}
\end{align*}
 where
 \begin{equation*}
 N_{j,p}^{k,W}(\rho,\phi) := L_q^{2jp+\gamma}\left(\frac{\rho^2}{2}\right)P_j^{l_1, l_0}(\cos(2p\phi)). 
 \end{equation*}
This suggests that the $W$-invariant generalized Hermite polynomials are given by
\begin{equation*}
H_{\tau_1,\tau_2}^W(\rho,\phi) = \sqrt{\frac{q!}{\Gamma(2jp+q+\gamma+1)}} \left(\frac{\rho^2}{2}\right)^{jp}N_{j,p}^{k,W}(\rho,\phi)
\end{equation*}
for $\tau_1 = 2q \,(q \geq 0), \tau_2 = 2jp \,(j \geq 0)$ and zero otherwise. An elegant proof of this claim was given to us (private communication) by Professor C. F. Dunkl and is as follows: the $j$-th $W$-invariant harmonic is given by (see Proposition 3. 15 in \cite{Dunkl0}) 
\begin{equation*}
h_j^W(\rho,\phi) = \rho^{2jp}P_j^{l_1, l_0}(\cos(2p\phi))
\end{equation*}
so that by Proposition 3.9 in \cite{Dunkl01}, 
\begin{equation*}
e^{-\Delta_k/2}[\rho^{2q}h_j^W(\rho,\phi)] = e^{-\Delta_k^W/2}[\rho^{2q}h_j^W(\rho,\phi)] = (-2)^j j! L_q^{2jp+\gamma}\left(\frac{\rho^2}{2}\right)P_j^{l_1, l_0}(\cos(2p\phi)).
\end{equation*}
\begin{nota}
A similar result holds for odd dihedral systems with $k_1=0, k_0 = k, p=n$.  
\end{nota}  

\section{Skew-product decomposition}
In this section, we derive a skew-product decomposition for $X$ using only independent Bessel processes. This is done by relating the process $\theta$ to a Jacobi process. That is why some results on Jacobi processes are collected below. 
\subsection{Facts on Jacobi processes} 
The Jacobi process $J$ of parameters $d,d' \geq 0$ is a $[0,1]$-valued process and is a solution, whenever it exists, of (\cite{War}) 
\begin{equation}\label{JSDE}
dJ_t = 2\sqrt{J_t(1-J_t)}dB_t + (d - (d+d')J_t) dt,
\end{equation}
where $B$ is a real Brownian motion.  As for squared Bessel processes, (\ref{JSDE}) has a unique strong solution for all $t \geq 0$ and all $J_0 \in [0,1]$ since the diffusion coefficient is $1/2$-H\"olderian and the drift term is Lipschitzien (\cite{Rev} p. 360). When $d \wedge d' \geq 2$ and $J_0 \in [0,1]$, then $J$ remains in $]0,1[$ while when 
$d\wedge d' > 0$, $J$ is valued in $[0,1]$ (\cite{Dou} p. 135). Besides, $J$ has the skew-product decomposition below (\cite{War}): 

\begin{equation}\label{WY}
\left(\frac{Z_1^2(t)}{Z_1^2(t) + Z_2^2(t)}\right)_{t \geq 0} = (J_{F_t})_{t \geq 0}, \quad F_t = \int_0^t \frac{ds}{Z_1^2(s)+Z_2^2(s)},
\end{equation}
where $Z_1,Z_2$ are two independent Bessel processes of dimension $d,d'$ respectively such that $ d+d' \geq 2$. Moreover, $J$ is independent from $Z_1^2+ Z_2^2$ (thereby from $F$). 

\subsection{Relating $\theta$ and $J$}  \label{S62}
Assume $d \wedge d' \geq 1$, and define 
\begin{equation*}
(H_t : = -\cos 2M_t)_{t \geq 0}
\end{equation*} 
where $(M_t = p\theta_{t/p^2})_{t \geq 0}$. Then an application of It\^o's formula and the pathwise uniqueness for the above SDE shows that $(H_t)_{t \geq 0}  = (Y_{2t})_{t \geq 0}$ where
\begin{equation*}
dY_t = \sqrt{2}\sqrt{1-Y_t^2}dB_t  - [(k_1 - k_0) + (k_0+k_1+1)Y_t]dt.
\end{equation*}
In fact, it is easy to see that $(1-Y_{2t})/2 = (1-H_t)/2 = \cos^2(M_t)$ is a Jacobi process of parameters $d = 2k_1+1,d' = 2k_0+1$. As a result, one gets
\begin{equation*}
\left(\theta_{A_t}\right)_{t \geq 0} = \left(\frac{1}{p}\arccos(\sqrt{J_{p^2A_t}})\right)_{t \geq 0}
\end{equation*}
where $J$ is independent from $X$, thereby from $A$.

\subsection{Skew-product decomposition} 
On the one hand, it is a well known fact (\cite{Rev}) that the sum of two independent squared Bessel processes of dimensions $d = 2k_1+1,d' = 2k_0+1$ is again a squared Bessel process of dimension $d+d'$, thus $Z := Z_1^2 + Z_2^2$ is a squared Bessel process of index $k_0+k_1$. On the other hand, for any conjuguate numbers $r,q$ and any Bessel process $R_{\nu}$ of index $\nu > -1/q$, there exists a Bessel process $R_{\nu q}$ of index $\nu q$ and defined on the same probability space such that the following holds (\cite{Rev})
\begin{equation}\label{time}
q^2R_{\nu}^{2/q}  = R_{\nu q}^2\left(\int_0^{\cdot}R_{\nu}^{-2/r}(s) ds\right).
\end{equation} 
Specializing (\ref{time}) with $\nu = k_0+k_1, q=p, R_{\nu} = Z$, there exists a Bessel process $R_{\nu q}$ such that 
\begin{equation}\label{time1}
Z^{2/p}_t = \frac{R_{\nu q}^2}{p^2}\left(\int_0^{t}\frac{ds}{Z_s^{2(p-1)/p}} \right) := \frac{1}{p^2} R_{\nu q}^2(\tau_t), \quad r = \frac{p}{p-1}.
\end{equation} 
Let $J$ be the Jacobi process defined in (\ref{WY}) with $d = 2k_1+1, d'=2k_0+1$ and define a radial Dunkl process $X$ by 
\begin{equation*}
X := \left(R_{\nu q}, \frac{1}{p}\arccos\sqrt{J}_{\displaystyle \int_0^{\cdot} p^2\frac{ds}{R_{\nu q}^2(s)}}\right) = \left(|X|, \theta_{A_t}\right)
\end{equation*}

Let $L_t := \inf \{s, \tau_s > t\}$ be the inverse of $\tau$ so that $Z^{2/p}_{L} = (1/p^2)|X|^2$. Then 
\begin{equation}\label{time2}
p^2A = p^2  \int_0^{\cdot} \frac{ds}{|X_s|^2} = \int_0^{\cdot} \frac{ds}{Z_{L_s}^{2/p}} = \int_0^{L_{\cdot}}\frac{d\tau_s}{Z_s^{2/p}}
= \int_0^{L_{\cdot}}\frac{ds}{Z_s^2} = F_{L_{\cdot}}.
\end{equation} 
As a result, when $k_0,k_1 \geq 0$ and $\theta \in C$, one has
\begin{equation*}
\left(\theta_{A_t}\right)_{t \geq 0}  = \left(\frac{1}{p}\arccos(\sqrt{J_{F_{L_t}}})\right)_{t \geq 0} =  \frac{1}{p}\left(\arccos\sqrt{\frac{Z_1^2}{Z_1^2 + Z_2^2}}(L_t)\right)_{t \geq 0}.
\end{equation*}
Finally 
\begin{pro}
\label{Pro1}
Let $k_0,k_1 \geq 0$ and define the time-change $\tau$ by
\begin{equation*}
\tau := \int_0^t \frac{ds}{[Z_1^2(s) + Z_2^2(s)]^{2(p-1)/p}}
\end{equation*}
where $Z_1, Z_2$ are two independent Bessel processes of dimensions $d=2k_1 + 1, d'=2k_0 +1$ respectively. 
Then there exists a radial Dunkl process $X$ associated with the even Dihedral group $\mathcal{D}_2(2p)$ such that $X_{\tau}$ is realized as the two-dimensional process 
\begin{equation*}
\left[p(Z_1^2 + Z_2^2)^{1/2p}, \frac{1}{p}\arccos\sqrt{\frac{Z_1^2}{Z_1^2 + Z_2^2}}\right].
\end{equation*}
A similar representation holds when $n$ is odd.
\end{pro}

\section{On the first hitting time of a wedge}
Let $X_0 = x \in C$ and let
\begin{equation*}
T_0 := \inf\{t, \, X_t \in \partial C\}
 \end{equation*} 
 be the first hitting time of $\partial C$. Recall that for dihedral groups $\mathcal{D}_2(2p)$, $C$ is a wedge of angle $\pi/(2p)$. Recall also from (\cite{Dem2}) that if the index function $l:= k-1/2$ takes one striclty negative value for some simple root $\alpha$, then $<\alpha,X>$ hits zero a.s. so that $T_0 < \infty$ a.s (see also \cite{Chy}). For even dihedral systems, two cases are to be considered: 
 \begin{itemize}
 \item $1/2 \leq k_0,k_1  \leq 1$ with either $k_1 > 1/2$ or $k_0 > 1/2$ or equivalently $0 \leq l_0,l_1 \leq 1/2$ with either $l_0 > 0$ or $l_1 > 0$: in that case, the radial Dunkl process with index function $-l$ hits $\partial C$ a.s. and we will use results from \cite{Chy}. 
 \item One and only one of the index values is strictly negative while the other value is positive: in that case, the radial Dunkl process of index function $l$ hits $\partial C$. We will follow  a different strategy based on our representation of the angular process $\theta$ as a Jacobi process and on results from (\cite{Dou}) on Jacobi processes. 
This strategy applies to the first case too.
\end{itemize}  
For odd dihedral systems, we can only have $1/2 < k \leq 1$ and computations are similar to the ones done in the first case for even dihedral systems.   

\subsection{Even dihedral groups: first case} 
We give two different approaches: while the first one has the advantage to be short, the second approach is shown to be efficient for both cases. In fact, the first approach uses the absolute-continuity relations for radial Dunkl processes derived in \cite{Chy} and we are met with a complicated exponential functional when dealing with the second case (see \cite{Dem} for more details). The second approach focus on the angular process $\theta$ which was identified with a Jacobi process and uses absolute-continuity relations for Jacobi processes from \cite{Dou}.
\\     
{\it First approach}: write $x =\rho e^{i\phi}, \rho > 0, 0 < \phi < \pi/(2p)$ and let $1/2 \leq k_0,k_1 \leq 1$ with either $k_0 > 1/2$ or $k_1 > 1/2$. 
Using part (c) of Proposition 2. 15 in \cite{Chy}, p. 38, the tail distribution of $T_0$ is given by: 
\begin{align*}
\p_x^{-l}(T_0 > t) = \e_x^l\left[\left(\prod_{\alpha \in R_+}\frac{<\alpha,X_t>}{<\alpha,x>}\right)^{-2l(\alpha)}\right],
\end{align*} 
where $\p_x^l$ ($\e_x^l$) denotes the law of a radial Dunkl process starting at $x \in C$ and of index $l$ (the corresponding expectation). From (\ref{RSG}), one gets 
\begin{equation*}
\p_x^{-l}(T_0 > t) = \frac{e^{-\rho^2/2t}}{c_k} \left(\frac{\rho}{\sqrt t}\right)^{2p(l_0+l_1)}\sin^{2l_0}(p\phi)\cos^{2l_1}(p\phi)g\left(\frac{\rho}{\sqrt t}, \phi \right),
\end{equation*}
where 
\begin{equation*}
g(\rho,\phi) = \int_0^{\infty}\int_0^{\pi/n}e^{-r^2/2}D_k^W(\rho,\phi,r,\theta) r^{2p+1}\sin(2p\theta)dr d\theta.
\end{equation*} 
With regard to Corollary \ref{GBF}, it amounts to evaluate
\begin{eqnarray*}
S_1(j)  &=&\int_0^{\infty}e^{-r^2/2} {\it I}_{b_j}(\rho r)r^{2p+1-\gamma}dr,\\
S_2 (j) &=& \int_0^{\pi/2p}P_j^{k_1-1/2,k_0-1/2}(\cos(2p\theta))\sin(2p\theta)d\theta
\end{eqnarray*}
for every $j \geq 0$, where $b_j := 2jp+\gamma$. In order to evaluate $S_1$, we use the expansion (\cite{Leb}, p. 108)
\begin{equation*}
{\it I}_{b_j}(\rho r) = \sum_{q \geq 0}\frac{1}{\Gamma(b_j+q+1)} \left(\frac{\rho r}{2}\right)^{2q+b_j}
\end{equation*}
and exchange the order of integration to get 
\begin{equation*}
S_1(j) = 2^{(p-\gamma)/2}\frac{\Gamma(a_j+1)}{\Gamma(b_j+1)} \left(\frac{\rho}{\sqrt 2}\right)^{b_j} {}_1\mathscr{F}_1\left(a_j+1,b_j+1, \frac{\rho^2}{2}\right)
\end{equation*}
where 
\begin{equation*}
a_j= \frac{(2j+k_0+k_1)p + 2p - \gamma}{2} = (j+1)p.
\end{equation*}
Using the variable change $ s = \cos(2p\theta)$, $S_2(j)$ transforms to 
\begin{align*}
S_2(j) = \frac{1}{2p}\int_{-1}^{1}P_j^{k_1-1/2,k_0-1/2}(s)ds = \frac{1}{p}\int_0^{1}P_j^{k_1-1/2, k_0-1/2}(2s-1)ds
\end{align*} 
which is easily computed using the expansion p. 21 in \cite{Dunkl}. As a result, the tail distribution is given by 
\begin{pro}\label{Pro2}
\begin{align*}
\p_x^{-l}(T_0 > t) &= \frac{e^{-\rho^2/2t}}{c_k} \left(\frac{\rho}{\sqrt t}\right)^{\gamma - 2p}\sin^{2l_0}(p\phi)\cos^{2l_1}(p\phi)
\\& \sum_{j \geq 0} S_2(j)\frac{\Gamma(a_j+1)}{\Gamma(b_j+1)} \left(\frac{\rho}{\sqrt{2t}}\right)^{2jp} {}_1\mathscr{F}_1\left(a_j+1,b_j+1, \frac{\rho^2}{2t}\right)P_j^{l_1, l_0}(\cos(2p\phi))
\\& = \frac{1}{c_k} \left(\frac{\rho}{\sqrt t}\right)^{\gamma -2p}\sin^{2l_0}(p\phi)\cos^{2l_1}(p\phi)
\\& \sum_{j \geq 0} S_2(j)\frac{\Gamma(a_j+1)}{\Gamma(b_j+1)} \left(\frac{\rho}{\sqrt{2t}}\right)^{2jp} {}_1\mathscr{F}_1\left(b_j - a_j,b_j+1, -\frac{\rho^2}{2t}\right)
P_j^{l_1, l_0}(\cos(2p\phi))
\end{align*}
by Kummer's transformation (\cite{Leb}). 
\end{pro}
\begin{nota}
The value $k \equiv 1$ corresponds to the first exit time of a Brownian motion from a wedge and our result fits the one in \cite{Ban}. Moreover $\gamma = 2p, b_j  =2 a_j$ and one may use the duplication formula to simplify the above ratio of Gamma functions and use some argument simplifications for the confluent hypergeometric function (\cite{Leb}). More precisely,  the Gauss duplication formula gives
\begin{equation*}
\sqrt{\pi} \Gamma(2a_j+1) = 2^{2a_j} \Gamma(a_j+1/2)\Gamma(a_j+1).
\end{equation*}
Moreover, 
\begin{equation*}
{}_1\mathscr{F}_1(a, 2a+1; z) = 2^{2a-1}\Gamma\left(a+ \frac{1}{2}\right)e^{z/2}(-z)^{1/2 - a}
\left[{\it I}_{a-1/2}\left(-\frac{z}{2}\right) + {\it I}_{a+1/2}\left(-\frac{z}{2}\right)\right]
\end{equation*}
and $P_j^{1, 1}(\cos(2p\phi))= (\sqrt{2}/\sqrt{\pi})\sin (2a_j \phi)/\sin(2p\phi)$ so that 
\begin{equation*}
S_2 (j) = (\sqrt{2}/\sqrt{\pi}) \int_0^{\pi/2p} \sin(2a_j\theta)d\theta = \frac{\sqrt{2}}{2a_j\sqrt{\pi}}(1- \cos[(j+1)\pi]) 
\end{equation*}
which equals zero when $j$ is odd and equals $\sqrt{2}/(\sqrt{\pi}a_j)$ when $j$ is even. Thus
\begin{equation*}
\p_x^{-1/2}(T_0 > t) = ce^{-\rho^2/4t} \left(\frac{\rho^2}{2t}\right)^{1/2-p}\sum_{j=0}^{\infty}\left[{\it I}_{d_j}\left(\frac{\rho^2}{4t}\right) + {\it I}_{d_j+1}\left(\frac{\rho^2}{4t}\right)\right]
\sin[2(2j+1)p\phi]
\end{equation*}
where $d_j = (2j+1)p-1/2$.
\end{nota}

\noindent {\it Second approach}: recall that 
\begin{equation*}
\left(\theta_{A_t}\right)_{t \geq 0} = \left(\frac{1}{p}\arccos(\sqrt{J_{p^2A_t}})\right)_{t \geq 0}
\end{equation*}
where $J$ is a Jacobi process of parameters $d = 2k_1 +1, d'= 2k_0 +1$ and is independent from $|X|$ (thereby from $A$). Then, for an appropriate index function $l$ 
(so that $T_0 < \infty$), one has 
\begin{align*}
\p_x^l(T_0 > t) &= \p_x^l(0 < \theta_{A_t} < \pi/(2p)) = \p_x^l(0 < J_{p^2A_t} < 1)  = \p_{x}^l(T_J > p^2 A_t)
\end{align*}
where 
\begin{equation*}
T_J := \inf\{t, \, J_t = 0\} \wedge \inf\{t, \, J_t = 1\}
\end{equation*}
is the first exit time from the interval $[0,1]$ by a Jacobi process. Now, let us recall the following absolute continuity relation between the laws of Jacobi processes of different set of parameters (Theorem 9.4.3. p.140 in \cite{Dou} specialized to $m=1$): let $\p_{z}^{d,d'}$ denote the probability law of a Jacobi process  starting at $z \in ]0,1[$ and of parameters $d \wedge d' > 0$. Writing $T$ for $T_J$, then 
\begin{align}\label{Theorem}
\p_{z}^{d_1,d_1'}{}_{|\mathscr{F}_t \cap \{T  > t\}} &= \left(\frac{J_t}{z}\right)^{\kappa}\left(\frac{1-J_t}{1-z}\right)^{\beta}
\exp -\int_0^t ds\left[c + \frac{u}{J_s} + \frac{v}{1-J_s}\right] \p_{z}^{d_2,d_2'}{}_{|\mathscr{F}_t \cap \{T  > t\}} 
\end{align}      
where $(\mathscr{F}_t)_t$ is the natural filtration of $J$, $d_i \wedge d_i'> 0, i=1,2$ and
\begin{eqnarray*}
\kappa &=& \frac{d_1-d_2}{4}, \beta =  \frac{d_1'-d_2'}{4}, \\
u &=& \frac{d_1-d_2}{4}\left(\frac{d_1+d_2}{2}-2\right), v = \frac{d_1'-d_2'}{4}\left(\frac{d_1'+d_2'}{2}-2\right),\\
c &=& \frac{d_1 + d_1' - d_2 - d_2'}{4}\left(2-  \frac{d_1 + d_1' +d_2 + d_2'}{2}\right). 
\end{eqnarray*} 
One Corollary of the above Theorem (Corollary 9.4.6. p.140\footnote{The exponential factor $e^{-ct}$ was missed in \cite{Dou}}.) states that if $d := 2k_1 + 1 := 2(l_1 +1), d' := 2k_0 +1:= 2(l_0 +1)$ with $0 \leq l_1,l_0 < 1$, then 
\begin{equation}\label{Relation1}
\p_{z}^{-l_1,-l_0}(T > t) = \e_z^{l_1,l_0}\left[e^{-ct}\left(\frac{z}{J_t}\right)^{l_1}\left(\frac{1-z}{1-J_t}\right)^{l_0}\right].
\end{equation}
where $\p_z^{l_1,l_0}$ ($\e_z^{l_1,l_0}$) denotes the probability law of a Jacobi process of index $l_1,l_0$ and starting at $z$ (the corresponding expectation).\\
To recover the result in Proposition \ref{Pro2}, proceed as follows. Let $-1/2 \leq -l_0,-l_1\leq 0$ so that at least one value is strictly negative. Using the semi-group density 
of $J = \cos^2M$ which follows from (\ref{Jacobi}): 
\begin{align*}
p_t^J(z,s) &= \frac{1}{2\sqrt{s(1-s)}}m_t^k(\arccos \sqrt{z}, \arccos \sqrt{s})\\&
= c_k\left[\sum_{j\geq 0}e^{\lambda_j t}P_j^{l_1, l_0}(2z-1)P_j^{l_1, l_0}(2s-1)\right] s^{l_1}(1-s)^{l_0}, \lambda_j = -2j(j+k_0+k_1),
\end{align*} 
for some constant $c_k$ and $z,s \in ]0,1[$, together with the independence of $J$ and $A$ and (\ref{Relation1}), one gets 
\begin{align*}
\p_{x}^{-l}&(T_0 > t) =
\e_z^{l_1,l_0}\left[e^{-cp^2A_t}\left(\frac{z}{J_{p^2A_t}}\right)^{l_1}\left(\frac{1-z}{1-J_{p^2A_t}}\right)^{l_0}\right]
\\&= c_kz^{l_1}(1-z)^{l_0}\sum_{j \geq 0}\e_{\rho}^{\gamma'}[e^{-p^2(c-\lambda_j)A_t}]P_j^{l_1, l_0}(2z-1)\int_0^1P_j^{l_1, l_0}(2s-1)ds,
 \end{align*}  
where $\e_{\rho}^{\gamma'}$ is the law of the Bessel process $|X|$ of index 
\begin{equation*}
\gamma' = (2-k_0-k_1)p = 2p - \gamma,
\end{equation*}
corresponding to $-l$ and $c = 2(l_0+l_1) = 2(k_0 + k_1 -1)$. Now, note that integral in the RHS is up to a constant $S_2(j)$ defined in the previous section and that $z^{l_1}(1-z)^{l_0} = \sin^{2l_0}(p\phi)\cos^{2l_1}(p\phi)$. Next, use the formula   
\begin{equation*}
\mathbb{E}_{\rho}^{\gamma'}[e^{-(c-\lambda_j) p^2A_t}| |X_t|=r]  = \frac{{\it I}_{\sqrt{\gamma'^2 + 2(c-\lambda_j)p^2}}(\rho r/t)}{{\it I}_{\gamma'}(\rho r/t)},
\end{equation*}
 and the semi group density of $|X|$ 
\begin{equation*}
q_t(\rho,r) = \frac{1}{t} \left(\frac{r}{\rho}\right)^{\gamma'}r e^{-(\rho^2+r^2)/2t} {\it I}_{\gamma'}\left(\frac{\rho r}{t}\right)
\end{equation*}
to deduce after elemantary computations that 
\begin{align*}
\mathbb{E}_{\rho}^{\gamma'}[e^{-(c-\lambda_j) p^2A_t}] &= \frac{e^{-\rho^2/2t}}{t} \rho^{2p-\gamma} \int_0^{\infty} r^{2p-\gamma+1}
e^{-r^2/2t} {\it I}_{2jp + \gamma}\left(\frac{\rho r}{t}\right)dr  
\\& =  e^{-\rho^2/2t} \left(\frac{\rho}{\sqrt{t}}\right)^{\gamma-2p} \int_0^{\infty} r^{2p-\gamma+1}e^{-r^2/2} {\it I}_{2jp + \gamma}\left(\rho r\right)dr
\\& =  e^{-\rho^2/2t} \left(\frac{\rho}{\sqrt{t}}\right)^{\gamma-2p}S_1(j) 
\end{align*}
In fact,  
\begin{align*}
\gamma'^2 - 2(c-\lambda_j)p^2 &= [(2-k_0-k_1)^2 + 4[(k_0+k_1-1) + j(j+k_0+k_1)]]p^2 \\&
= [(k_0+k_1)^2 - 4(k_0+k_1-1) +4(k_0+k_1-1)  +4j(j+k_0+k_1)]p^2 
\\&= (2j+k_0+k_1)^2p^2 = (2jp+\gamma)^2.
\end{align*}
Finally, it only remains to relate the modified Bessel function ${\it I}_{\nu}$ and the hypergeometric function ${}_0\mathscr{F}_1$ via: 
\begin{equation*}
{\it I}_{\nu}(z) = \frac{1}{\Gamma(\nu+1)}z^{\nu}{}_0\mathscr{F}_1(\nu+1,z).
\end{equation*}

\subsection{Even dihedral groups: second case}
We use the second approach developed above and we suppose for instance that $k_1 < 1/2$ while $k_0 \geq 1/2$. Take $1 \leq d_1 = d < 2, d_1'=d' \geq 2$ in (\ref{Theorem})
then perform the parameters change  
\begin{equation*}
d_2' = d_1' = d' \geq 2, \, d_2 = 4- d_1 = 4- d = 3-2k_1 > 2 
\end{equation*}   
so that the indices corresponding to the new parameters $d_2,d_2'$ are positive for which $T = T_J = \infty$ a.s.. Moreover, one has $\beta = u = v = 0$ which yields 
\begin{align*}
\p_{x}^l(T_0> t) &= \p_z^{4-d, d'} \left[e^{-cp^2A_t}\left(\frac{J_{p^2A_t}}{z}\right)^{\kappa}\right]
\end{align*}  
where $\kappa = (d- d_2)/2 = d/2-1= l_1 < 0$ and $c = -d'(d-2) = d'(2-d) > 0$. Since the parameter $d_2$ corresponds to the index value $1/2-k_1 = -l_1> 0$ and multiplicity function $1-k_1 > 1/2$, we get
\begin{equation*}
\p_{x}^l(T_0>t) = c_kz^{1-2k_1}(1-z)^{l_0}\sum_{j \geq 0}\e_{\rho}^{\gamma}[e^{p^2(\lambda_j-c)A_t}]P_j^{1/2-k_1, k_0-1/2}(2z-1)F(j) 
\end{equation*}
where $\lambda_j = -2j(j+k_0 + 1-k_1), \gamma = (k_0+k_1)p$ and   
\begin{equation*}
F(j) = \int_0^1(1-s)^{l_0} P_j^{1/2-k_1, k_0-1/2}(2s-1) ds.
\end{equation*}
We leave the computations to the interested reader. A similar result holds when $k_0 <1/2$ and $k_1 \geq 1/2$: substitute $k_1,l_1$ by $k_0,l_0$ respectively and $s, z$ by $1-s, 1-z$ respectively.
\\  

{\bf Acknowedgments:} the author would like to give a special thank for Professor C. F. Dunkl for his intensive reading of the paper, for his fruitful remarks and for pointing to him  the references \cite{Dunkl0}, \cite{Dunkl01}. The author also thanks Professor M. Bozejko for his remarks and for stimulating discussions at the Wroclaw Institute of Mathematics.


\begin{thebibliography}{99}
\bibitem{Baez}\emph{J. C. Baez}. The octonions. {\it Bull. Amer. Math. Soc. (N. S.)} {\bf 39}, no. 2. 2002, 145-205.
\bibitem{Bak}\emph{T. H. Baker, P. J. Forrester}. The Calogero-Sutherland model and generalized classical polynomials. {\it Comm. Math. Phys}. {\bf 188}. 1997, 175-216.
\bibitem{Bakry}\emph{D. Bakry}. Remarques sur les semi-groupes de Jacobi. {\it Hommage \`a P. Andr\'e Meyer et J. Neveu. Ast\'erisque}. {\bf 236}, 1996, 23-39.
\bibitem{Ban}\emph{R. Ba$\tilde{n}$uelos, R. G. Smits}. Brownian motions in cones. {\it P. T. R. F. }{\bf 108}. 1997, 299-319.
\bibitem{Chy}\emph{O. Chybiryakov}. Processus de Dunkl et Relation de Lamperti. Ph. D. Thesis, Paris VI Univ. June 2005.
\bibitem{De}\emph{N. Demni, M. Zani}. Large deviations for statistics of Jacobi process. {\it To appear in S. P. A.}
\bibitem{Dem}\emph{N. Demni}. First hitting time of the boundary of a Weyl chamber by radial Dunkl processes. {\it SIGMA Journal.} {\bf 4}, 2008, 074, 14 pages.
\bibitem{Dem1}\emph{N. Demni}. Generalized Bessel function of type $D$.  {\it SIGMA Journal.} {\bf 4}, 2008, 075, 7 pages.
\bibitem{Dem2}\emph{N. Demni}. Note on radial Dunkl processes. {\it Submitted to Ann. I. H. P.}
\bibitem{Dou}\emph{Y. Doumerc}. Matrix Jacobi Process. {\it Ph. D. Thesis}. Paul Sabatier Univ. May 2005. 
\bibitem{Dunkl}\emph{C. F. Dunkl, Y. Xu}. Orthogonal Polynomials of Several Variables. {\it Encyclopedia of Mathematics and Its Applications. Cambridge University Press}. 2001.
\bibitem{Dunkl0}\emph{C. F. Dunkl}. Differential-difference operators associated to reflection groups. {\it Trans. Amer. Math. Soc.} {\bf 311}. 1989, no. 1, 167-183.
\bibitem{Dunkl01}\emph{C. F. Dunkl}. Integral kernels with reflection group invariance. {\it Canad. J. Math.} {\bf 43}. 1991, no. 6, 1213--1227.
\bibitem{Dunkl1}\emph{C. F. Dunkl}. Generating functions associated with Dihedral groups. {\it Special functions (Hong Kong 1999), World Sci. Publ. Rier Edge, NJ}. 2000, 72-87.
\bibitem{Gra}\emph{D. J. Grabiner}. Brownian motion in a Weyl chamber, non-colliding particles and random matrices. {\it Ann. IHP}. {\bf 35}, 1999, no. 2. 177-204. 
\bibitem{Hum}\emph{J. E. Humphreys}. Reflections Groups and Coxeter Groups. {\it Cambridge University Press}. {\bf 29}. 2000.
\bibitem{Leb}\emph{N. N. Lebedev}. Special Functions And Their Applications. {\it Dover Publications, INC}. 1972.
\bibitem{Rev}\emph{D. Revuz, M. Yor}. Continuous Martingales And Brownian Motion, $3^{\textrm{rd}}$ ed, Springer, 1999.
\bibitem{Schou}\emph{W. Schoutens}. Stochastic Processes and Orthogonal Polynomials. {\it Lecture Notes in Statistics}, {\bf 146}. Springer, 2000.
\bibitem{Ros}\emph{M. R\"osler}. Dunkl operator : theory and applications, orthogonal polynomials and special functions (Leuven, 2002). {\it Lecture Notes in Math}. Vol. 1817, Springer, Berlin, 2003, 93-135.
\bibitem{War}\emph{J. Warren, M. Yor}. The Brownian Burglar : conditioning Brownian motion by its local time process. {\it S\'em. Probab.} {\bf XXXII.}, 1998, 328-342.
\bibitem{Yor}\emph{M.Yor}. Loi de l'indice du lacet brownien et distribution de Hartman-Watson, {\it Zeit. Wahr. Verw. Geb.} {\bf 53}, no.1, 1980, 71-95.
\end{thebibliography}
\end{document}